\documentclass[a4paper,10pt]{article}
\usepackage{amsmath,amsfonts,amssymb,amsthm} 

\newcommand{\Rset}{\mathbb{R}}
\newcommand{\Nset}{\mathbb{N}}

\newcommand{\Cset}{\mathbb{C}}
\newcommand{\Qset}{\mathbb{Q}}
\newcommand{\PP}{\mathbb{P}}
\newcommand{\EE}{\mathbb{E}}
\newcommand{\DS}{\displaystyle}

\newtheorem{Theorem}{Theorem}[section]
\newtheorem{Definition}[Theorem]{Definition}
\newtheorem{Proposition}[Theorem]{Proposition}
\newtheorem{Lemma}[Theorem]{Lemma}
\newtheorem{Corollary}[Theorem]{Corollary}
\newtheorem{Remark}[Theorem]{Remark}
\newtheorem{Example}[Theorem]{Example}

\newtheorem{Hypothesis}[Theorem]{Hypothesis}

\title{Kolmogorov equations for measures}
\author{Luigi Manca\\
Scuola Normale Superiore\\
Piazza dei Cavalieri, 7\\
56126, Pisa, Italy}
\begin{document}

\maketitle
\begin{abstract}
We consider a semigroup of operators in the Banach space $C_b(H)$ of uniformly continuous and bounded functions on a separable Hilbert space $H$.
In particular, we deal with semigroups that are related to solution of stochastic PDEs in $H$ and which are not, in general, strongly continuous.
We prove an existence and uniqueness result for a measure valued equation involving this class of semigroups.
Then we apply the result to a large class of second order  differential operators in $C_b(H)$.
\end{abstract}

\section{Introduction}
Let $H$ be a separable Hilbert space (with norm $|\cdot|$ and inner product $\langle\cdot,\cdot\rangle $), and let $\mathcal B(H)$ be its Borel $\sigma$-algebra.
We are concerned with semigroups of operators on $C_b(H)$, the Banach space of all uniformly continuous and bounded functions $f:H\to\Rset$, endowed with the supremum norm $\|\cdot\|_0$. 
In particular, we consider a semigroup of linear and bounded operators $\{P_t\}_{t\geq0} \subset \mathcal L(C_b(H) )$ which is a {\em stochastically continuous Markov  semigroup},
that is there exists a family $\{\pi_t(x,\cdot),\,t\geq0,\,x\in H\}$ of probability Borel measures on $H$ such that
\begin{itemize}
\item the map $\Rset^+\times H\to [0,1]$, $(t,x)\mapsto\pi_t(x,\Gamma)$ is measurable, for any Borel set $\Gamma\in \mathcal B(H)$; %
\item  $\DS \pi_{t+s}(x,\Gamma )=
\int_H\pi_s(y,\Gamma )\pi_t(x,dy)$,
 for all $t,s\geq 0, \; x \in H,\,\Gamma  \in \mathcal B(H)$;
\item  for any $x \in H$, $\pi_0(x,\cdot)=\delta_x(\cdot)$, the probability measure concentrated in $x$; 
\item $\DS P_t\varphi(x)=\int_H\varphi(y)\pi_t(x,dy) $, for any $t\geq0$, $\varphi \in C_b(H)$, $x\in H$;
\item for any  $\varphi \in C_b(H)$, $x\in H$, the function $\Rset^+\to\Rset$, $t\mapsto   P_t\varphi(x)$ is continuous. 
\end{itemize}

It is well known that such semigroups are related with the solution of stochastic partial differential equations in $H$, see \cite{DP04}, \cite{DPZ92}, \cite{DPZ02},  \cite{Manca06}.
It is also well known that, in general, they are not strongly continuous in $C_b(H)$ (see, for instance, \cite{Cerrai}, \cite{GK01}, \cite{Priola}).
However, we can define an infinitesimal generator $(K,D(K))$ by setting 
\begin{equation}  \label{e.3}
\begin{cases}
  \DS D(K)=\bigg\{ \varphi \in C_b(H): \exists g\in C_b(H), \lim_{t\to 0^+} \frac{P_t\varphi(x)-\varphi(x)}{t}=g(x),  \\ 
  \DS  \qquad\qquad\quad x\in H,\;\sup_{t\in(0,1)}\left\|\frac{P_t\varphi-\varphi}{t}\right\|_0<\infty \bigg\}\\   
   {}   \\
  \DS K\varphi(x)=\lim_{t\to 0^+} \frac{P_t\varphi(x)-\varphi(x)}{t},\quad \varphi\in D(K),\,x\in H.
\end{cases}
\end{equation}
%
%
%
%
Denoting by $\mathcal M(H)$ the space of all finite Borel measures on $H$, we consider the following problem for measures: 
{\em given $\mu\in  \mathcal M(H)$, find a family of Borel finite measures $\{\mu_t\}_{t\geq0}$ such that }
\begin{equation}\label{e.2.5}
\begin{cases} \DS
\frac{d}{dt} \int_H\varphi(x)\mu_t(dx) = \int_H K  \varphi(x)\mu_t(dx) & t\geq0,\,\varphi \in D(K )\\
\DS\mu_0 = \mu,& \mu\in {\cal M}(H).
\end{cases}
\end{equation}
To give a precise meaning of this problem, we introduce the notion of solution of \eqref{e.2.5}
\begin{Definition}\label{d.1.1}{\em
Given  $\mu \in {\cal M}(H)$, we say that a family of measures $\{\mu_t\}_{t\geq 0}$ is a {\em solution of the  measure equation} \eqref{e.2.5} if   the following is fulfilled 
\begin{itemize}
\item the total variation of the measures $\mu_t$ satisfies
\begin{equation} \label{e.5b}
   \int_0^T\|\mu_t\|_{TV}dt<\infty,\quad T>0; 
\end{equation} 
\item for any $\varphi \in D(K)$, the real valued function 
\begin{equation} \label{e.fefe}
   \Rset^+\to \Rset,\qquad t\mapsto \int_H \varphi(x)\mu_t(dx)
\end{equation}
is absolutely continuous\footnote{that is it belongs to $W^{1,1}([0,T])$} and for any $t\geq0$  it holds 
\begin{equation} \label{e.2.6}
 \int_H \varphi(x)\mu_t(dx)-\int_H \varphi(x)\mu (dx)=\int_0^t\bigg(\int_H K\varphi(x)\mu_s(dx)\bigg)ds. 
\end{equation} 
\end{itemize}
}
\end{Definition}

The first result of this paper is the following
\begin{Theorem} \label{t.1}
Let $\{P_t\}_{t\geq0}$ be a stochastically continuous Markov semigroup and let $(K,D(K))$ be its infinitesimal generator, defined as in \eqref{e.3}.
Then, the formula 
$$
  \langle \varphi, P_t^*F\rangle_{\mathcal L(C_b(H),\,(C_b(H))^*)} = \langle P_t\varphi, F\rangle_{\mathcal L(C_b(H),\, (C_b(H))^*)}
$$
defines a semigroup $(P_t^*)_{t\geq0}$ of linear and continuous operators  on $(C_b(H))^*$ that maps $\mathcal M(H)$ into $\mathcal M(H)$.
Moreover, for any $\mu\in \mathcal M(H)$, $\varphi\in C_b(H)$ the map 
\begin{equation}  \label{e.fifi}
   \Rset^+\to \Rset,\quad t\mapsto \int_H\varphi(x) P_t^*\mu(dx)
\end{equation}
is continuous, and if $\varphi\in D(K)$ it is also differentiable with continuous differential
\begin{equation}  \label{e.fifi2}
 \frac{d}{dt} \int_H \varphi(x)P_t^*\mu(dx) = 
   \int_H K \varphi(x)P_t^*\mu(dx).
\end{equation}
Finally, for any $\mu\in  \mathcal M(H)$ there exists a unique solution of   the measure equation \eqref{e.2.5}, given by $\{P_t^*\mu\}_{t\geq0}$.  
\end{Theorem}
%
%


%
%
%
In second part of this paper, we consider the transition semigroup $\{P_t\}_{t\geq0}$ associated to the  stochastic differential equation in $H$ 
\begin{equation} \label{e.3.1}
\left\{\begin{array}{lll}
dX(t)&=& \big(AX(t)+F(X(t))\big)dt+Q^{1/2}dW(t),\quad t\geq 0\\
\\
X(0)&=&x\in H,
\end{array}\right.
\end{equation}
where 
\begin{Hypothesis}  \label{h.3.1}
\begin{itemize}

\item[(i)] $A\colon D(A)\subset H\to H$ is the infinitesimal
generator of a strongly continuous semigroup $e^{tA}$ of type ${\cal G}(M,\omega)$, i.e. there exist $M\geq 0$ and $\omega\in \Rset$ such that $\|e^{tA}\|_{{\mathcal L}(H)}\le Me^{\omega t}$, $ t\ge 0$;
\item[(ii)]  $Q\in \mathcal L(H)$ is non negative and symmetric, so its square root $Q^{1/2}$ exists and it is unique (cf, for instance, \cite{Schmeidler}). 
Moreover, for any $t> 0$ the linear operator $Q_{t}$, defined by  
\begin{equation*}  
Q_{t}x=\int_{0}^{t}e^{sA}Qe^{sA^{*}}xds,\;\;x\in H,\;t\ge 0
\end{equation*}
has finite trace;
\item[(iii)] $F:H\to H$ is a Lipschitz continuous map;
\item[(iv)] $(W(t))_{t\geq0}$ is a cylindrical Wiener process, defined on a stochastic basis $(\Omega,\mathcal F, (\mathcal F_t)_{t\geq0},\PP)$ and with values in $H$.
\end{itemize}
\end{Hypothesis}
It is well known that under Hypothesis \ref{h.3.1} equation \eqref{e.3.1} has  a unique stochastically continuous  {\em mild} solution $(X(t,x))_{t\geq0,x\in H}$ (see, for instance, \cite{DPZ92}), that is the random variable $X(t,x):\Omega\to H$ is solution of the integral equation
\begin{equation} \label{e.4.3}
 X(t,x)= e^{tA}x+\int_0^te^{(t-s)A}Q^{1/2}dW(t)+\int_0^te^{(t-s)A}F(X(s,x))ds
\end{equation}
and that 
\begin{equation} \label{e.10}
    \lim_{t\to t_0}\EE\big[|X(t,x)-X(t_0,x)|^2\big]=0, 
\end{equation}
for any $t_0 \geq0 $.
Hence,  the transition semigroup $\{P_t\}_{t\geq0}$ in $C_b(H)$ associated to equation \eqref{e.3.1} is defined by setting
\begin{equation} \label{e.38a}
   P_t\varphi(x)=\EE\big[\varphi(X(t,x))\big],\quad \varphi \in C_b(H),\, t\geq 0,\, x\in H.
\end{equation}
It is not too hard to prove that $\{P_t\}_{t\geq0}$ is a stochastically continuous Markov semigroup (cf Proposition \ref{p.4.2}).
 This allows us to define the infinitesimal generator $(K,D(K))$ of $\{P_t\}_{t\geq0}$, as in \eqref{e.3}.
We are interested in the relationships between  $(K,D(K))$ and the {\em Kolmogorov} differential operator
\begin{equation} \label{e.4}
 K_0\varphi(x)= \frac12\textrm{Tr}\big[QD^2\varphi(x)\big]+\langle x,A^*D\varphi(x)\rangle+\langle D\varphi(x),F(x)\rangle,\,x\in H.
\end{equation}
In order to study this problem, we shall introduce the notions of $\pi$-convergence and of $\pi$-core (cf section 2).
Then we shall prove the second result of this paper
\begin{Theorem} \label{t.2} 
Let $\mathcal I_A(H)$ be the linear span of the real and imaginary part of the functions
$$
H\to \Cset,\quad x\mapsto \int_{0}^{a}e^{i\langle
e^{sA}x,h \rangle -\frac12 \langle Q_sh,h\rangle}ds:\; a>0,\;h\in D(A^{*}), 
$$
where $D(A^*)$ is the domain of the adjoint operator of $A$.
Then $ \mathcal I_A(H)\subset D(K)$ and for any $\varphi\in \mathcal I_A(H)$  we have $K\varphi=K_0\varphi$.
Moreover, the set $ \mathcal I_A(H)$ is a $\pi$-core for $(K,D(K))$.
\end{Theorem}
The theorem above states, in particular, that $(K,D(K))$ is an extension of $K_0$.
The problem of extending a differential operator of the form \eqref{e.4} to an infinitesimal generator of a diffusion semigroup has been the object of many papers
in the recent years. 
For instance, when the semigroup has an invariant measure $\nu$, this problem can be studied in the Hilbert space $L^2(H;\nu)$ of all  Borel function $f:H\to\Rset$ which are square integrable with respect to $\nu$ (see, for instance, \cite{DPZ96}, \cite{DPZ02}, \cite{DP04}, \cite{Manca06} and references therein).
Other similar results have been stated  by studying $\{P_t\}_{t\geq0}$ in weighted spaces (see, for instance, \cite{DPT01},  \cite{DP04} and references therein).

Results about this problem in $C_b(H)$ are, at our knowledge, new.
As consequence of Theorem \ref{t.2} we have the third main result
\begin{Theorem} \label{t.1.4}
For any $\mu\in\mathcal M(H)$ there exists an unique solution $\{\mu_t\}_{t\geq0}\subset \mathcal M(H)$ of the measure equation \begin{equation}\label{e.2.5a}
\begin{cases} \DS
\frac{d}{dt} \int_H\varphi(x)\mu_t(dx) = \int_H K_0  \varphi(x)\mu_t(dx) & t\geq0,\,\varphi \in\mathcal I_A(H)\\
\DS\mu_0 = \mu,& \mu\in {\cal M}(H)
\end{cases}
\end{equation}
and this solution  is done by $\mu_t=P_t^*\mu$.
\end{Theorem}

Kolmogorov equations for measures have been the object of several papers.
Recently, by starting with a generalization of the classical work of Hasminskii (see the monograph \cite{HAS80}),
in \cite{BR00} has been stated sufficient conditions in order to ensure existence of 
a weak solution for 
partial differential operators of the form 
$$
    H\varphi(t,x)= a^{ij}(t,x)\partial_{x_i}\partial_{x_j}\varphi(x)+b^i(t,x)\partial_{x_i}\varphi(x),\quad (t,x)\in (0,1)\times \Rset^d,
$$
where $\varphi\in C_0^\infty( R^d) $ and $a^{ij},b^i\colon (0,1)\times R^d\to \Rset$ are suitable locally integrable functions.
The authors proves that if there exists a Lyapunov-type function for the operator $H$, then there exists a probability measure on $\Rset^d$ that solves the equation $H^*\nu=0$\index{}, that is
$$
   \int_{\Rset^d}H\varphi(x)\nu(dx)=0
$$
for any test function $\varphi\in C_0^\infty( R^d)$.
In \cite{BR01}, this result has been extended to separable Hilbert spaces.
With similar techniques, in  \cite{BDPR04} the results have been extended to parabolic differential operators of the form $Lu(t,x)=u_t(t,x)+Hu(t,x)$,
$u\in  C_0^\infty((0,1)\times \Rset^d)$.
The authors proved that if there exists a Lyapunov-type function for the operator $L$, then for any probability measure $ \nu$ on $\Rset^d$ there exists a family of probability measures $\{\mu_t\}_{t\in(0,1)}$ such that
$$
    \int_0^1\int_{R^d} Lu(t,x)\mu_t(dx)dt=0
$$ 
for any $u\in C_0^\infty((0,1)\times \Rset^d)$ and $\lim_{t\to0}\int_{\Rset^d} \zeta(x)\mu_t(dx)=\int_{\Rset^d} \zeta(x)\nu(dx)$, for any $\zeta\in C_0^\infty( R^d)$.

In our paper, we are concentrated in uniqueness of the solution.
Indeed, we deal with differential operators that are related to diffusion processes,  hence it is not difficult to prove existence of a solution.
To get uniqueness we need, of course, suitable regularity properties of the coefficients.

Uniqueness results for such a kind of differential operators in Hilbert spaces  are, at our knowledge, new.
In a forthcoming paper we shall study the case of reaction-diffusion, Burgers and Navier-Stokes operators.

Let us describe how is organize this paper.
In the next section we introduce notations and prove some  results about approximation of $C_b(H)$ functions by trigonometric series and some properties of the solutions of the measure equation \eqref{e.2.5}.
In section 3, we prove Theorem \ref{t.1}. 
Finally, Theorems \ref{t.2} and \ref{t.1.4} are proved in section 4 and 5, respectively.  
In order to be  clear, each proof is divided into several subsections.

\section{Notations and preliminary results}
We shall identify $H$ with its topological dual space $H^*$.
If $E$ is a Banach space, we denote by $C_b(H;E)$ the Banach space of all uniformly continuous and bounded functions $f:H\to E$, endowed the supremum norm $\|\cdot\|_{C_b(H;E)}$.  
$\mathcal L(H;E)$ is the usual Banach space of all the linear and continuous operators  $A:H\to E$, endowed with the norm $\|\cdot\|_{\mathcal L(H;E)}$.
If $E=\Rset$, we briefly write $\mathcal L(H)$ instead of $\mathcal L(H;\Rset)$.
$C_b^1(H;E)$ denotes the space of all the functions $f\in C_b(H;E)$ which are Fr\'echet differentiable with uniformly continuous and bounded differential $DF\in C_b(H;\mathcal L(H;E))$.
As above, we shall use the notation $C_b^1(H)=  C_b^1(H;\Rset)$.
Inductively, $C_b^k(H;E)$ is the Banach space of all functions $f\in C_b^{k-1}(H;E)$ which are $k$-times Fr\'echet differentiable with uniformly continuous and bounded differential.

We deal with semigroup of operators that are not strongly continuous.
For this reason, we introduce the notion of $\pi$-convergence in the space $C_b(H)$ (see \cite{Priola}).
\begin{Definition}{\em
A sequence $\{\varphi_n\}_{n\in \Nset} \subset C_b(H)$ is said to be {\em $\pi$-convergent} to a function $\varphi$ $\in$ $C_b(H)$  if for any $x\in H$ we have 
$$
 \lim_{n\to \infty}\varphi_n(x)=\varphi(x)
$$ 
and 
$$
  \sup_{n\in\Nset}\|\varphi_n\|_0 <\infty. 
$$
Similarly,  the $m$-indexed sequence $\{\varphi_{n_1,\ldots,n_m}\}_{n_1\in\Nset,\ldots,n_m\in\Nset}\subset C_b(H)$ is said to be $\pi$-convergent to  $\varphi$ $\in$ $C_b(H)$ if for any $i\in \{2,\ldots,m\}$  
 there exists an $i-1$-indexed sequence $\{\varphi_{n_1,\ldots,n_{i-1}}\}_{n_1\in\Nset,\ldots,n_{i-1}\in\Nset}\subset C_b(H)$  such that
$$
 \lim_{n_i\to \infty}\cdots\lim_{n_m\to \infty} \varphi_{n_1,\ldots,n_m}(x)=\varphi_{n_1,\ldots,n_{i-1}}(x),\quad x\in H
$$ 
and 
$$
 \lim_{n_i\to\infty}\varphi_{n_1,\ldots,n_i}\stackrel{\pi}{=}\varphi_{n_1,\ldots,n_{i-1}}.
$$
We shall write
$$
 \lim_{n_1\to \infty}\cdots\lim_{n_m\to \infty} \varphi_{n_1,\ldots,n_m}\stackrel{\pi}{=}\varphi
$$ 
or $\varphi_n\stackrel{\pi}{\to}\varphi$ as $n \to \infty$, when the sequence has one index. 
}
\end{Definition}
\begin{Remark}{\em
Note that  the $\pi$-convergence implies the convergence in $L^p(H;\mu)$, for any $\mu\in{\cal M}(H)$, $p\in [1,\infty)$.}
\end{Remark}
\begin{Remark}{\em
The notion of $\pi$-convergence is considered also in \cite{EK}, under the name of {\em boundedly and pointwise} convergence.
}
\end{Remark}
\begin{Remark}{\em
The topology on $C_b(H)$ induced by the $\pi$-convergence is not sequentially complete.
For a survey on this fact see \cite{GK01}, \cite{Priola} .
}
\end{Remark}
\begin{Definition}{\em
For any subset $D\subset C_b(H)$ we say that $\varphi$ belongs to the $\pi$-closure of $D$, and we  denote it by $\varphi \in \overline{D}^\pi$, if there exists $m\in\Nset$ and an $m$-indexed sequence $\{\varphi_{n_1,\ldots,n_m}\}_{n_1\in\Nset,\ldots,n_m\in\Nset}\subset D$ such that
$$
   \lim_{n_1\to\infty}\cdots\lim_{n_m\to\infty}\varphi_{n_1,\ldots,n_m}\stackrel{\pi}{=}\varphi.
$$
Finally, we shall say that a subset $D\subset C_b(H)$ is $\pi$-dense in $C\subset C_b(H)$ if $\overline{D}^\pi = C$.
}
\end{Definition}
\begin{Remark}{\em
In order to avoid heavy notations, we shall often assume  that the sequences have only one index.
}
\end{Remark}
It will be helpful the following results about approximation of $C_b(H)$-functions by trigonometric series.
\begin{Proposition} \label{p.2.1}
We denote by $\mathcal E (H) $  the linear span  of the real and imaginary part of the functions
$$
 H\to\Cset,\quad    x\mapsto e^{i \langle x, h\rangle },
$$ 
where $h\in H$.
Then $\mathcal E (H) $ is $\pi$-dense in $C_b(H)$ and for any $\varphi\in C_b(H)$ there exists a two-indexed sequence $(\varphi_{n_1, n_2})\subset \mathcal E(H) $ such that
\begin{eqnarray} 
&\DS \lim_{n_1\to \infty }\lim_{n_2\to \infty }\varphi_{n_1, n_2}(x)=\varphi(x),\, x\in H & \label{e.2.21}\\
&\DS\sup_{n_1, n_2}\| \varphi_{n_1, n_2}\|_0\leq \|\varphi\|_0.&  \label{e.2.22}
\end{eqnarray}
Moreover, if $\varphi\in C_b^1(H)$ we can choose the sequence $(\varphi_{n_1, n_2})\subset \mathcal E(H) $ in such a way that \eqref{e.2.21}, \eqref{e.2.22} hold and for any $h\in H$ 
\begin{eqnarray} 
&\DS \lim_{n_1\to \infty }\lim_{n_2\to \infty }\langle D\varphi_{n_1, n_2}(x),h\rangle=\langle D\varphi(x),h\rangle,\, x\in H & \notag\\
&\DS\sup_{n_1, n_2}\|D \varphi_{n_1, n_2}\|_{C_b(H;H)}\leq \|D\varphi\|_{C_b(H;H)}.&  \label{e.2.24}
\end{eqnarray}
\end{Proposition}
\begin{proof}
In  \cite[Proposition 1.2]{DP04} are proved \eqref{e.2.21}, \eqref{e.2.22}.
\eqref{e.2.24} follows by the well known properties of the Fourier approximation with F\'ejer kernels of differentiable functions (see, for instance, \cite{KO}).
\end{proof}

\begin{Proposition} \label{p.2.4}
Let $(e_k)_{k\in\Nset}$ be a complete orthonormal system of $H$.
We denote by $\mathcal E_\Qset(H) $  the $\Qset$-linear span  of the real and imaginary part of the functions
$$
     x\mapsto e^{i \langle x, q_1e_1+\cdots+ q_n e_n\rangle },
$$ 
where $n\in \Nset$ and $q_1,\cdots,q_n \in\Qset$.
Then, $\mathcal E_\Qset(H) $ is a countable $\pi$-dense subset of $C_b(H)$.
In particular, for any $\varphi\in C_b(H)$ there exists a three-indexed sequence $(\varphi_{n_1, n_2, n_3})\subset \mathcal E(H) $ such that
\begin{eqnarray} 
\DS \lim_{n_1\to \infty }\lim_{n_2\to \infty }\lim_{n_3\to \infty }\varphi_{n_1, n_2, n_3}\stackrel{\pi}{=}\varphi. &\label{e.21}  
\end{eqnarray}
\end{Proposition}
\begin{proof}
By Proposition \ref{p.2.1} we can find a two-indexed sequence $(\varphi_{n_1, n_2})\subset \mathcal E(H)$ such that \eqref{e.2.21}, \eqref{e.2.22} holds.
It is now clear that  we can approximate any $\varphi_{n_1, n_2}$ by 
 a three-indexed sequence $(\varphi_{n_1, n_2, n_3})\subset \mathcal E_\Qset(H) $ such that $\forall n_1,n_2\in\Nset$ it holds  
\begin{eqnarray*} 
&\DS   \lim_{n_3\to \infty }\varphi_{n_1, n_2, n_3}(x)=\varphi_{n_1, n_2}(x),\, x\in H&\\
& \DS  \sup_{ n_3\in\Nset}\|\varphi_{n_1, n_2, n_3}\|_0 < \infty.& .
\end{eqnarray*}
Hence, it follows that the three-indexed sequence $(\varphi_{n_1, n_2, n_3})$ is the claimed one.
\end{proof}

It turns out that  a stochastically continuous Markov semigroup  
 is a $\pi$-semigroup, as introduced  by Priola (see \cite{Priola}).
So, we have the following
\begin{Proposition} \label{p.1}
Let $\{P_t\}_{t\geq0}$ be a stochastically continuous Markov semigroup and let $(K,D(K))$ be its infinitesimal generator, defined as in \eqref{e.3}.
Then,  $\{P_t\}_{t\geq0}$  is  a $\pi$-semigroup, that is
\begin{itemize}
\item[(i)]  for any $t\geq0$, $P_t\in \mathcal L (C_b(H))$ and $\|P_t\|_{{\cal L}(C_b(H))}\leq 1$;
\item[(ii)] $P_tP_s=P_{t+s}$, $t,s \geq 0$;
\item[(iii)] for any $\varphi\in C_b(H)$ and $x\in H$, the map $[0,\infty)\to\Rset$, $t\mapsto P_t\varphi(x)$ is continuous;
\item[(iv)] for any sequence $\{\varphi_n\}_{n\in\Nset} \subset C_b(H)$ such that $\varphi_n\stackrel{\pi}{\to}\varphi$ as $n\to\infty$ we have $P_t\varphi_n\stackrel{\pi}{\to}P_t\varphi$, as $n\to\infty$, for any $t\geq0$.
\end{itemize}
\end{Proposition}

The next result is proved in \cite[Propositions 3.2, 3.3, 3.4]{Priola}.
\begin{Theorem} \label{t.2.9}
Under the assumption of Proposition \ref{p.1},  we have
\begin{itemize}
\item[(i)] for any $\varphi \in D(K)$, $P_t\varphi \in D(K)$ and $KP_t\varphi = P_tK\varphi$, $t\geq0$;
\item[(ii)] for any $\varphi \in D(K)$, $x\in H$, the map $[0,\infty)\to \Rset$, $t\mapsto P_t\varphi(x)$ is continuously differentiable and $(d/dt)P_t\varphi(x) = P_tK\varphi(x)$;
\item[(iii)] $D(K)$ is $\pi$-dense in $C_b(H)$;
\item[(iv)] $K$ is a $\pi$-closed operator on $C_b(H)$, that is for any sequence $\{\varphi_n\}_{n\in\Nset}\subset C_b(H)$ such that  $\varphi_n\stackrel{\pi}{\to}\varphi\in C_b(H)$ and $K\varphi_n\stackrel{\pi}{\to}g\in C_b(H)$ as $n\to \infty$ it follows that $\varphi \in D(K)$ and $g=K\varphi$;
\item[(v)] for any $f\in C_b(H)$, $t>0$ the map 
$
      H\to \Rset,\quad x\mapsto \int_0^t P_sf(x) ds
$
    belongs to $D(K)$ and it holds
\begin{equation*} 
K\left( \int_0^tP_sf ds\right) = P_tf-f.
\end{equation*}
Moreover, if $\varphi\in D(K)$ we have
$$
K\left( \int_0^tP_sf ds\right) =  \int_0^tKP_sf ds;
$$
\item[(vi)] for any $\lambda>0$ the linear operator $R(\lambda,K)$ on $C_b(H)$ done by
$$
   R(\lambda,K)f(x)=\int_0^\infty e^{-\lambda t}P_tf(x)dt, \quad f\in C_b(H),\,x\in H
$$
satisfies, for any  $f\in C_b(H)$
$$
   R(\lambda,K)\in \mathcal L(C_b(H)),\quad \quad \|R(\lambda,K)\|_{\mathcal L(C_b(H))}\leq \frac{1}{\lambda} 
$$
$$
  R(\lambda,K)f\in D(K),\quad (\lambda I-K)R(\lambda,K)f=f.
$$
We call $R(\lambda,K)$ the {\em resolvent} of $K$ at $\lambda$.
\end{itemize}
\end{Theorem}
%
%
A central role will be play by the notion of {\em $\pi$-core}. 
\begin{Definition}
We shall say that a set $D\subset D(K)$ is a {\em $\pi$-core} for the operator $(K,D(K))$ if $D$ is $\pi$-dense in $C_b(H)$ and for any $\varphi\in D(K)$ there exists $m\in\Nset$ and an  $m$-indexed sequence $\{\varphi_{n_1,\ldots,n_m}\}_{n_1\in\Nset,\ldots,n_m\in\Nset} \subset D$ such that
$$
 \lim_{n_1\to \infty}\cdots\lim_{n_m\to \infty} \varphi_{n_1,\ldots,n_m}\stackrel{\pi}{=}\varphi
$$ 
and
$$
 \lim_{n_1\to \infty}\cdots\lim_{n_m\to \infty} K\varphi_{n_1,\ldots,n_m}\stackrel{\pi}{=}K\varphi.
$$
\end{Definition}
It is clear that a  $\pi$-core in nothing but the extension of the notion of {\em core} with respect to the $\pi$-convergence.
An useful example of core is done by the following 
\begin{Proposition} \label{p.4.3}
 Let $ \{P_t\}_{t\geq0}$ be a  stochastically continuous Markov semigroup and let $(K,D(K))$ be its infinitesimal generator.
If $D\subset D(K)$ in $\pi$-dense in $C_b(H)$ and $P_t (D)\subset D$ for all $t\geq0$, then $D$ is a $\pi$-core for $(K,D(K))$.
\end{Proposition}
\begin{proof}
In order to get the result, we proceed as in \cite{EN00}.
Let  $\varphi \in D(K)$.
Since $D$ in $\pi$-dense in $C_b(H)$, there exists a  sequence $(\varphi_{n_2})\subset D$ (for the sack of simplicity we assume that the sequence has only one index) such that $\varphi_{n_2}\stackrel{\pi}{\to}\varphi$ as $n_2\to\infty$. 
Set
\begin{equation} \label{e.32a}
   \varphi_{n_1,n_2,n_3}(x)=\frac{1}{n_3}\sum_{i=1}^{n_3}P_{\frac{i}{n_1n_3}}\varphi_{n_2}(x)
\end{equation}
for any $n_1,n_2,n_3\in \Nset$.
By Hypothesis, $(\varphi_{n_1,n_2,n_3})\subset D$.
Taking into account Proposition \ref{p.1},  a strightforward  computation shows  that for any $x\in H$
\begin{eqnarray*}
   \lim_{n_1\to\infty}\lim_{n_2\to\infty}\lim_{n_3\to\infty}\varphi_{n_1,n_2,n_3}(x)
 &=& \lim_{n_1\to\infty}\lim_{n_2\to\infty}  n_1\int_0^{\frac{1}{n_1}}P_t\varphi_{n_2}(x)dt\\
 &=&  \lim_{n_1\to\infty}n_1\int_0^{\frac{1}{n_1}}P_t\varphi(x)dt=\varphi(x).
\end{eqnarray*}
Moreover,
$$
    \sup_{n_1,n_2,n_3\in\Nset} \|\varphi_{n_1,n_2,n_3}\|_0\leq \sup_{n_2}\|\varphi_{n_2}\|_0<\infty
$$
since $\varphi_{n_2}\stackrel{\pi}{\to}\varphi$ as $n_2\to\infty$.
Hence, 
$$
  \lim_{n_1\to\infty}\lim_{n_2\to\infty}\lim_{n_3\to\infty}\varphi_{n_1,n_2,n_3}\stackrel{\pi}{=}\varphi.
$$
Similarly, since $D\subset D(K)$ and Theorem \ref{t.2.9} holds, we have
\begin{eqnarray*}
   \lim_{n_3\to\infty}K\varphi_{n_1,n_2,n_3}(x) &=& n_1\int_0^{\frac{1}{n_1}}KP_t\varphi_{n_2}(x)dt \\
  &=&n_1\Big(P_{\frac{1}{n_1}}\varphi_{n_2}(x)-\varphi_{n_2}(x)\Big).
\end{eqnarray*}
So we find
$$
  \lim_{n_1\to\infty}\lim_{n_2\to\infty}\lim_{n_3\to\infty}K\varphi_{n_1,n_2,n_3}(x)
   =\lim_{n_1\to\infty}\lim_{n_2\to\infty}n_1\Big(P_{\frac{1}{n_1}}\varphi_{n_2}(x)-\varphi_{n_2}(x)\Big)
$$
\begin{equation} \label{e.33}
  = \lim_{n_1\to\infty}n_1\Big(P_{\frac{1}{n_1}}\varphi(x)-\varphi(x)\Big)=K\varphi(x),
\end{equation}
since $\varphi\in D(K)$.
To conclude the proof, we have to show that these limits are uniformly bounded with respect to every index.
Indeed we have
$$
   \sup_{n_3\in\Nset}\|K\varphi_{n_1,n_2,n_3}\|\leq \|K\varphi_{n_2}\|<\infty,
$$
$$
    \sup_{n_2\in\Nset}\|n_1\Big(P_{\frac{1}{n_1}}\varphi_{n_2}-\varphi_{n_2}\Big)\|_0\leq 
     2 n_1  \sup_{n_2\in\Nset}\|\varphi_{n_2}\|_0<\infty.
$$
Finally, the last limit in \eqref{e.33} is uniformly bounded with respect to $n_1$ since $\varphi\in D(K)$.
\end{proof}
%
%
%
%
%
%
\section{Proof of Theorem \ref{t.1}} 
We split the proof in several lemma, collected into three parts:
in the first one we prove the first statement of the theorem;
in the second one, we prove existence of a solution;
finally, in the third part, we prove uniqueness of the solution.

\subsection{$P_t^*:\mathcal M(H)\to \mathcal M(H)$}

\begin{Lemma} \label{l.1}
Let $\{P_t\}_{t\geq0}$ be a  stochastically continuous Markov semigroup.
The family of  linear maps $\{P_t^*\}_{t\geq0}:(C_b(H))^*\to (C_b(H))^*$, defined by the formula
\begin{equation} \label{e.5}
  \langle \varphi, P_t^*F\rangle_{\mathcal L(C_b(H),\,(C_b(H))^*)} = \langle P_t\varphi, F\rangle_{\mathcal L(C_b(H),\, (C_b(H))^*)},\, 
\end{equation}
where $t\geq0, \, F\in (C_b(H))^*,\, \varphi \in C_b(H)$, is a  semigroup of linear maps  on $(C_b(H))^*$ of norm $1$ and maps $\mathcal M(H)$ into $\mathcal M(H)$.
\end{Lemma}
\begin{proof} 
Clearly, $P_t^*$ is linear.
Let $F\in \big(C_b(H)\big)^*$, $t\geq0$.
We have,  for any $\varphi \in C_b(H)$,
$$
   \langle \varphi, P_t^*F\rangle_{\mathcal L(C_b(H),\,(C_b(H))^*)} \leq \|\varphi\|_0\|F\|_{(C_b(H))^*}.
$$
Then $P_t:(C_b(H))^*\to (C_b(H))^*$ has norm equal to $1$.
Moreover, by \eqref{e.5} it follows easily that $P_t^*(P_s^*F)=P_{t+s}^*F$, for any $t,s\geq0$, $F\in (C_b(H))^*$.
Hence, \eqref{e.5} defines a semigroups of application in  $(C_b(H))^*$ of norm equal to $1$.

Now we prove that $P_t^*:\mathcal M(H)\to \mathcal M(H)$.
To check this, let $\{\pi_t(x,\cdot),\,x\in H\}$ be the family of probability measures associated to $P_t$, that is
$    P_t\varphi(x)=\int_H\varphi(y) \pi_t(x,dy)$, 
for any $\varphi \in C_b(H)$.
Hence, if $\mu \in \mathcal M(H)$, we can define the map $\Lambda:\mathcal B(H)\to [0,\infty)$ by
\begin{equation} \label{e.6}
   \Lambda(\Gamma)= \int_H \pi_t(x,\Gamma)\mu(dx), \quad \Gamma \in \mathcal B(H).
\end{equation}
It is easy to see that $\Lambda$ is a $\sigma$-additive Borel finite measure on $H$.
In order to conclude the proof, we shall show that $\langle \varphi, P_t^*\mu\rangle=\int_H\varphi(x)\Lambda(dx)$, for all $\varphi\in C_b(H)$. 
To see this, we extend the operator $P_t$ to an linear and continuous operator in $L^\infty(H;\Rset)$, still denoted by $P_t$.
This extension  follows by approximating pointwise any function $f\in L^\infty(H;\Rset)$ by a sequence of functions in $C_b(H)$.
Moreover, this extension is unique.
In particular, if $\Gamma$ is a Borel set of $H$, we have $P_t\chi_\Gamma(x)=\pi_t(x,\Gamma)$, $\forall x\in H$.
So, if $\varphi\in C_b(H)$, we can find a sequence $(\varphi_n)\subset L^\infty(H;\Rset)$ of functions of the form\footnote{$\chi_A$ is the characteristic function of the set $A\subset H$}
$$
   \varphi_n(x)=\sum_{k=1}^{N(n)}\alpha_k^n\chi_{A_k^n}(x),
$$ 
where $N(n)\in \Nset$,  $\alpha_k^n\in \Rset$, $A_k^n\in \mathcal B(H)$ are such that $A_k^n\cap A_l^n $ if $k\not= l$, $\bigcup_k A_k^n= H $ and
\begin{eqnarray*}
 \DS \sup_{n\in\Nset}\|\varphi_n\|_0\leq\|\varphi\|_0;\quad \lim_{n\to\infty}\varphi_n(x)=\varphi(x),\,x\in H
\end{eqnarray*}
Consequently, for  any $x\in H$ we have $P_t\varphi_n(x)=P_t\varphi(x)$ as $n\to\infty $ and $\sup_n\|P_t\varphi_n\|_0\leq \sup_n\| \varphi_n\|_0 \|\varphi\|_0$. 
By the dominated convergence theorem it follows
$$
   \int_H\varphi(x)\Lambda(dx)= \lim_{n\to\infty}\int_H\varphi_n(x)\Lambda(dx)=
$$
$$
 =\lim_{n\to\infty}\sum_{k=1}^{N(n)}\alpha_k^n\int_H\pi_t(x,A_k^n)\mu(dx)=   \lim_{n\to\infty}\int_HP_t\varphi_n(x)\mu(dx)=\langle \varphi,P_t^*\mu\rangle.
$$
Hence, the result follows.
\end{proof}

\subsection{Existence of a solution}
\begin{Lemma}
Let $\mu\in\mathcal M(H)$.
Under the hypothesis of Lemma \ref{l.1}, the semigroup $\{P_t^*\mu\}_{t\geq0}$ defined in \eqref{e.5} is a solution of the measure equation \eqref{e.2.5}.
Moreover, if $\varphi\in C_b(H)$ the function \eqref{e.fifi} is continuous, and if $\varphi\in D(K)$ function \eqref{e.fifi} is also differentiable with continuous differential done by \eqref{e.fifi2}.
\end{Lemma}
By Lemma \ref{l.1}, for any $\mu\in \mathcal M(H)$ the formula \eqref{e.5} define a family $\{P_t^*\mu\}_{t\geq0}$ of measures  on $H$.
Since  for any $\varphi\in C_b(H)$ it holds
$$
   \lim_{t\to 0^+} \int_HP_t\varphi(x)\mu(dx)=\int_H\varphi(x)\mu(dx),
$$ 
 by the semigroup property of $P_t$ it follows that for any  $\varphi \in C_b(H)$ the function
\begin{equation}  \label{e.7}
 \Rset^+ \to \Rset,\quad t\mapsto \int_H\varphi(x)P_t^*\mu(dx)
\end{equation}
is continuous.
Clearly, $P_0^*\mu=\mu$.
Now we show that if $\varphi \in D(K)$ then the function \eqref{e.7} is differentiable.
Indeed, by taking into account \eqref{e.3} and that $P_t^*\mu\in \mathcal M(H)$, for any  $\varphi \in D(K)$  we can apply the dominated convergence theorem to obtain 
\begin{eqnarray*}
  &&\frac{d}{dt}\int_H \varphi(x)P_t^*\mu(dx)=\\
  &&\qquad = \lim_{h\to 0}\frac1h \left(\int_H P_{t+h}\varphi(x)\mu(dx)- \int_HP_t\varphi(x)\mu_t(dx)\right)  \\
&&\qquad=\lim_{h\to 0}\int_H \left(\frac{P_{t+h}\varphi(x) - P_t\varphi(x)}{h}\right)\mu(dx)\\
 &&\qquad =\lim_{h\to 0}\int_H P_t\left(\frac{P_h\varphi - \varphi}{h}\right)(x)\mu(dx)\\
 &&\qquad =\int_H \lim_{h\to 0}\left(\frac{P_h\varphi - \varphi}{h}\right)(x)P_t^*\mu(dx)= 
  \int_H K\varphi(x)P_t^*\mu(dx).
\end{eqnarray*}
Then, by arguing as above, the differential of \eqref{e.7} is continuous. 
This clearly implies that   $\{P_t^*\mu\}_{t\geq0}$ is a solution of the measure equation \eqref{e.2.5}.
                                        \subsection{Uniqueness of the solution}
Since problem \eqref{e.2.5} is linear, it is enough to take $\mu = 0$. 
We claim that $\mu_t = 0$, $\forall t\geq0$.
In order to prove this, let us fix $T>0$ and let us consider the Kolmogorov backward equation
\begin{equation} \label{e.2.7}
\begin{cases}
u_t(t,x) + Ku(t,x)=\varphi(x) & t \in[0,T],\, x\in H, \\
u(T,x)=0, & {}
\end{cases}
\end{equation}
where  $\varphi\in C_b(H)$.  
The meaning of \eqref{e.2.7} is make clear by the following lemma. 
\begin{Lemma} \label{l.2.12}
 For any $T>0$,  $\varphi\in C_b(H)$ 
  the real valued function 
\begin{eqnarray} 
  &&u:[0,T]\times H\to \Rset \notag \\
  &&u(t,x)= -\int_0^{T-t}P_s\varphi(x)ds,\quad (t,x)\in [0,T]\times H.\label{e.2.8}
\end{eqnarray}
satisfies the following statements
 \begin{itemize}
  \item[{\em (i)}] $u\in C_b([0,T]\times H)$ \footnote{Clearly, $C_b([0,T]\times H)$ is isomorphic to $C([0,T]; C_b(H))$};
  \item[{\em (ii)}] $u(t,\cdot)\in D(K)$ for any $t\in[0,T]$ and the function $[0,T]\times H\to \Rset$, $(t,x)\mapsto Ku(t,x)$ is continuous and bounded; 
  \item[{\em (iii)}]  the real valued function $[0,T]\times H\to\Rset$, $(t,x)\mapsto u(t,x)$ is derivable with respect to $t$ with continuous and bounded derivative $u_t(t,x)$, that is for any $x$ the function $u(\cdot,x)$ is differentiable with differential $u_t(\cdot,x)$, and the function $[0,T]\times H\to \Rset$, $(t,x)\mapsto u_t(t,x)$ is continuous and bounded;
  \item[{\em (iv)}] for any $(t,x)\in [0,T]\times H$ the function $u$ satisfies  \eqref{e.2.7}.  
 \end{itemize}
\end{Lemma}
\begin{proof}
For any $s,t\in [0,T]$, $s\leq t$ we have
\begin{eqnarray*}
  u(t,x)-u(s,x)&=&-\int_0^{T-t}P_\tau\varphi(x)d\tau+\int_0^{T-t}P_\tau\varphi(x)d\tau\\
   &=& \int_{T-t}^{T-s}P_\tau\varphi(x)d\tau.
\end{eqnarray*}
Then
$$
  \|u(t,\cdot)-u(s,\cdot)\|_0\leq |t-s|\|\varphi\|_0.
$$
(i) is proved.
By (vi) of Theorem \ref{t.2.9}, $u(t,\cdot)\in D(K)$ for any $t\in [0,T]$ and it holds
$Ku(t,x)=-P_{T-t}\varphi(x)+\varphi(x)$, for any $x\in H$.
So  (ii) follows  (cf (iii) of Proposition \ref{p.1}).
Now let $h\in (-t,T-t)$ and $x\in H$. 
We have
\begin{eqnarray}
  &&\frac{u(t+h,x)-u(t,x)}{h}+Ku(t,x)-\varphi(x)=\label{e.8}\\
  &&\quad=\frac1h\int_{T-t-h}^{T-t}P_{T-s}\varphi(x)ds-P_{T-t}\varphi(x)\notag\\
  &&\quad =\frac1h\int_{T-t-h}^{T-t}\big(P_{T-s}\varphi(x)-P_{T-t}\varphi(x)\big)ds.\notag
\end{eqnarray}
Then, since $P_t\varphi(x)$ is continuous in $t$, \eqref{e.8} vanishes as $h\to0$.
This implies that $u(t,x)$ is derivable with respect to $t$ and \eqref{e.2.7} holds. 
Moreover, by (ii), we have that the maps $t\mapsto u_t(t,x)= -Ku(t,x)+\varphi(x)$ is continuous.
This proves (iii) and (iv).
The proof is complete.
\end{proof}

We need the following
\begin{Lemma} \label{l.2.13}
 Let $\{\mu_t\}$ be a solution of the measure equation \eqref{e.2.5} in the sense of Definition \ref{d.1.1}. 
 Then, for any function $u:[0,T]\times H\to \Rset $ satisfying statements (i), (ii), (iii) of Lemma \ref{l.2.12} the map
$$
   [0,T]\to \Rset,\qquad t\mapsto \int_Hu(t,x)\mu_t(dx)
$$
is absolutely continuous and for any $t\geq0$ it holds 
\begin{multline} \label{e.2.9}
  \int_Hu(t,x)\mu_t(dx)-\int_Hu(0,x)\mu(dx) \\
     =\int_0^t\left(\int_H\big(u_s(s,x)+Ku(s,x)\big)\mu_s(dx)\right)ds.
\end{multline}
\end{Lemma}
\begin{proof}
 We split the proof in several steps.\\
{\em Step 1: Approximation of $u(t,x)$.}\\
With no loss of generality, we   assume $T=1$.
For any $x\in H$, let us consider the approximating functions $\{u^n(\cdot,x)\}_{n\in\Nset}$ of 
$u(\cdot,x)$ done by the Bernstein polynomials (see, for instance, \cite{YO}). 
Namely, for any $n\in\Nset$, $x\in H$ we consider the function
$$
 [0,T]\to \Rset,\quad t\mapsto u^n(t,x)=\sum_{k=0}^n\alpha_{k,n}(t)u\Big(\frac{k}{n},x\Big),
$$
where 
$$
   \alpha_{k,n}(t)=\binom{n}{k}t^k(1-t)^{n-k}.
$$
Since $u\in C([0,T];C_b(H))$, it is well known that it holds 
\begin{equation} \label{e.2.10}
   \lim_{n\to\infty}\sup_{t\in[0,1]}\|u^n(t,\cdot)-u(t,\cdot)\|_0=0
\end{equation}
and
$$
   \sup_{t\in[0,1]}\|u^n(t,\cdot)\|_0<\infty,\quad n\in\Nset.
$$
Then, for any $t\in[0,1]$
\begin{equation}\label{e.2.11}
  \lim_{n\to\infty} u^n(t,\cdot)\stackrel{\pi}{=}u(t,\cdot). 
\end{equation}
We also have that for any $n\in\Nset$, $t\in[0,1]$
$$
   u^n(t,\cdot)\in D(K),
$$
and that for any $x\in H$ the function $[0,1]\to \Rset$, $t\mapsto Ku^n(t,x)$ is continuous (cf (ii) of Lemma \ref{l.2.12}).
Then, for any $x\in H$ it holds
$$
  \lim_{n\to\infty} \sup_{t\in [0,1]}|Ku^n(t,x)-Ku(t,x)|=0,
$$
\begin{equation}\label{e.25}
 \sup_{t\in[0,1]}\|Ku^n(t,\cdot)\|_0\leq \sup_{t\in[0,1]}\|Ku(t,\cdot)\|_0<\infty.
\end{equation}
This clearly implies that for any $t\in[0,1]$
\begin{equation}\label{e.2.12}
\lim_{n\to\infty}Ku^n(t,\cdot)\stackrel{\pi}{=}Ku(t,\cdot).
\end{equation}
Similarly, since for any $x$ the function $t\mapsto u(t,x)$ is differentiable with respect to $t$, we also have that for any $x\in H$
\begin{equation*} 
\lim_{n\to\infty} \sup_{t\in [0,1]}|u_t^n(t,x)-u_t(t,x)|=0,
\end{equation*}
\begin{equation}\label{e.27}
 \sup_{t\in[0,1]}\|u_t^n(t,\cdot)\|_0\leq \sup_{t\in[0,1]}\|u_t(t,\cdot)\|_0<\infty.
\end{equation}
Hence, for any $t\in [0,1]$ 
\begin{equation} \label{e.2.13}
 \lim_{n\to\infty}u_t^n(t,\cdot)\stackrel{\pi}{=}u_t(t,\cdot).
\end{equation}

\noindent
{\em Step 2: differential of $\int_Hu^n(t,x)\mu_t(dx)$}\\
For any  $n\in \Nset$, $k\leq n$ and for almost all $t\in[0,1]$ we have
\begin{eqnarray*}
   &&  \frac{d}{dt}\left(\int_H\alpha_{k,n}(t)u\Big(\frac{k}{n},x\Big)\mu_t(dx)\right)=\\
  &&\quad = \frac{d}{dt}\left(\alpha_{k,n}(t)\int_Hu\Big(\frac{k}{n},x\Big)\mu_t(dx)\right)\\
  &&\quad = \alpha_{k,n}'(t)\int_Hu\Big(\frac{k}{n},x\Big)\mu_t(dx)+
       \alpha_{k,n}(t)\int_HKu\Big(\frac{k}{n},x\Big)\mu_t(dx). \\
  &&\quad =\int_H\bigg(\alpha_{k,n}'(t)u\Big(\frac{k}{n},x\Big)+\alpha_{k,n}(t)Ku\Big(\frac{k}{n},x\Big)\bigg)\mu_t(dx).
\end{eqnarray*}
Note that the last terms belong to $L^1([0,1])$.
This implies
\begin{multline*}  
  \int_Hu^n(t,x)\mu_t(dx)-\int_Hu^n(0,x)\mu(dx) \\
     =\int_0^t\left(\int_H\big(u_s^n(s,x)+Ku^n(s,x)\big)\mu_s(dx)\right)ds,
\end{multline*}
for any $n\in\Nset$.

\noindent
            {\em Step 3: Conclusion}\\
Consider the functions
$$
   f:[0,1]\to \Rset, \quad f(t)=\int_Hu(t,x)\mu_t(dx)
$$
and  
$$
  f_n:[0,1]\to \Rset, \quad f_n(t)=\int_Hu^n(t,x)\mu_t(dx).
$$
By \ref{e.2.10}  we have
$$
    \left|\int_H\big(u^n(t,x)-u(t,x)\big)\mu_t(dx)    \right| 
   \leq  \sup_{t\in [0,1]}\| u^n(t,\cdot)-u(t,\cdot)   \|_0  \|\mu_t\|_{TV}. 
$$
Since \eqref{e.5b} and \eqref{e.2.10} hold, it follows that the sequence $(f_n)$ converges  to $f$ in $L^1([0,1])$, as $n\to\infty$.
We also have, by Step 2, that $f_n$ is absolutely continuous and hence differentiable in almost all $t\in[0,1]$, with differential in $L^1([0,1])$ done by
$$
 f_n'(t)= \int_H\big(u_t^n(t,x)+Ku^n(t,x)\big)\mu_t(dx),
$$
for almost all $t\in[0,1]$.
By \eqref{e.2.12}, \eqref{e.2.13} we have 
\begin{eqnarray} 
\lim_{n\to\infty}f_n'(t) &=&  \lim_{n\to\infty}\int_H\big(u_t^n(t,x)+Ku^n(t,x)\big)\mu_t(dx) \notag\\
 &=& \int_H\big(u_t (t,x)+Ku (t,x)\big)\mu_t(dx), \label{e.18}
\end{eqnarray}
for all $t\in [0,T]$.
Moreover, it holds
\begin{equation*}  
 \sup_{n\in\Nset} |f_n'(t)|\leq \bigg(\sup_{t\in [0,1]}\| u (t,\cdot)\|_0 +\sup_{t\in [0,1]}\|Ku (t,\cdot)\|\bigg) \|\mu_t\|_{TV}.
\end{equation*}
Hence, still  by \eqref{e.25}, \eqref{e.27}, there exists a constant $c>0$ such that $\sup_n |f_n'(t)| \leq c\|\mu_t\|_{TV}$.  
By taking into account \eqref{e.5b}, it follows that the limit in \eqref{e.18} holds in $L^1([0,1])$.
Let us denote by $g(t)$ the right-hand side of \eqref{e.18}.
We find, for any $a,b\in [0,1]$, 
\begin{eqnarray*}
    f(b)-f(a)
  &=& \lim_{n\to\infty}\big(f_n(b)-f_n(a)\big)\\
  &=& \lim_{n\to\infty}\int_a^bf'_n(t)dt
  = \int_a^b\lim_{n\to\infty}f'_n(t)dt=\int_a^bg(t)dt.
\end{eqnarray*}
Therefore, $f$ is absolutely continuous, and  $f'(t)=g(t)$ for almost all $t\in[0,1]$.
Lemma \ref{l.2.13} is proved.
\end{proof}
Now let $\varphi\in C_b(H)$ 
 and $u$ be the function defined in \eqref{e.2.8}.
We have that $u$ satisfies statements (i)--(iv) of Lemma \ref{l.2.12}. 
Hence, by Lemma \ref{l.2.13} it follows that the function $[0,T]\to \Rset$, $t\to \int_Hu(t,x)\mu_t(dx)$ is absolutely continuous, with differential 
\begin{eqnarray*}
  \frac{d}{dt}\int_H u(t,x) \mu_t(dx) &=& 
  \int_H\big( u_t(t,x)+Ku(t,x)\big) \mu_t(dx)\\
 &=&  \int_H\varphi(x)\mu_t(dx),  
\end{eqnarray*}
for almost all $t\in[0,T]$.
So, we can write
\begin{eqnarray*}
   0&=&\int_Hu(T,x)\mu_T(dx)-\int_Hu(0,x)\mu(dx)=\\
   &=& \int_0^T\left(\frac{d}{dt}\int_Hu(t,x)\mu_t(dx)\right)dt\\
    &=& \int_0^T\left(\int_H\varphi(x)\mu_t(dx)\right)dt.
\end{eqnarray*}
for all $\varphi \in C_b(H)$.
By the arbitrariness of $T$, it follows that for any $0\leq s\leq t$
it holds 
$$
     \int_s^t\left(\int_H\varphi(x)\mu_\tau(dx)\right)d\tau=0.
$$
Since \eqref{e.5b} holds, the function $t\mapsto \int_H\varphi(x)\mu_\tau(dx)$ belongs to $L^1([0,T])$, for any $T>0$. 
Consequently, by the well known properties of the Lebesgue integrable functions, for any $\varphi\in C_b(H)$ we have
\begin{equation} \label{e.20}
   \int_H\varphi(x)\mu_t(dx)=0,
\end{equation}
for almost all $t\geq 0$.
At this point, it is not clear if $\mu_t=0$ for almost all $t\geq 0$.
So, let us consider the set $\mathcal E_\Qset(H) $ introduced in Proposition \ref{p.2.4}.
We denote by $I_\varphi$ the set $\{ t\geq 0:\text{ \eqref{e.20} does}$ $\text{not hold} \}$  and by $I$ the set
$$
 I= \bigcup_{\varphi \in \mathcal E_\Qset(H) }I_\varphi.
$$
Since $\mathcal E_\Qset(H) $ is countable and for any $\varphi \in \mathcal E_\Qset(H) $ the set $I_\varphi$ is Borel and of Lebesgue measure equal to zero, then $I$ is Borel and of Lebesgue measure equal to zero.
It is clear that \eqref{e.20} holds for all $\varphi \in \mathcal E_\Qset(H) $, $t  \in \Rset^+\setminus I$. 
Now let $\varphi \in C_b(H)$. 
Still by Proposition \ref{p.2.4} we know that there exists a three-indexed sequence $(\varphi_{n_1, n_2, n_3})\subset \mathcal E_\Qset(H) $ such that \eqref{e.21} holds.
Hence, for any $t\in\Rset^+\setminus I$ we have
$$
  \int_H\varphi(x)\mu_t(dx)=\lim_{n_1\to \infty }\lim_{n_2\to \infty }\lim_{n_3\to \infty } \int_H\varphi_{n_1, n_2, n_3}(x)\mu_t(dx)=0.
$$
This implies that $\mu_t=0$ for all $t\in\Rset^+\setminus I$ and hence $\mu_t=0$ for almost all $t\geq0$.
The proof is now complete. \qed
\begin{Remark}{\em
In the last part of the proof it has a fundamental role the fact that the space $C_b(H)$ has a $\pi$-dense countable subset.
This is possible since $H$ is separable, as it can be see by Proposition \ref{p.2.4}.
} 
\end{Remark}

%


\section{Proof of Theorems \ref{t.1}}
We begin by showing that the transition semigroup $\{P_t\}_{t\geq0}$in \eqref{e.38a} is a stochastically continuous Markov semigroup in $C_b(H)$.
\begin{Proposition} \label{p.4.2}
Under  Hypothesis \ref{h.3.1}, the transition semigroup $\{P_t\}_{t\geq0}$ defined in \eqref{e.38a} is a  stochastically continuous Markov semigroup in $C_b(H)$.
\end{Proposition}
\begin{proof}
The fact that $\{P_t\}_{t\geq0}$ maps $C_b(H)$ into $C_b(H)$ and that it is a semigroup of operators  may be found in \cite[Proposition 3.9]{DP04}.
We also have $ P_t\varphi(x)= \int_H\varphi(y)\pi_t(x,dy)$, where  $\pi_t(x,\cdot)$ is the probability Borel measure on $H$ defined by $\pi_t(x,\Gamma)= \PP(X(t,x)\in \Gamma)$, $\forall \Gamma\in \mathcal B(H)$.
Hence, the semigroup $\{P_t\}_{t\geq0}$ is Markovian.
Finally, since $X(t,x)$ fulfills \eqref{e.10}, 
it follows easily that for any $\varphi\in C_bb(H) $,   $x\in H$ the function $H\to \Rset $, $t\to P_t\varphi(x)$   is  continuous.
\end{proof}

In order to prove Theorem \ref{t.1}, namely  $K\varphi=K_0\varphi$ if $\varphi\in \mathcal I_A(H)$ and that $\mathcal I_A(H)$ is a $\pi$-core for $(K,D(K))$,
we proceed by several steps. 
We start by studying the case when $F=0$ in \eqref{e.4}.

%
%
%
\subsection{The Ornstein-Uhlenbeck operator}
If $F=0$,  the operator \eqref{e.4} is known as the {\em Ornstein-Uhlenbeck} (OU) operator.
Let us consider the OU semigroup  $\{R_t\}_{t\geq0}$ done by  
\begin{equation*} 
 R_{t}\varphi (x)=\int_{H}\varphi (e^{tA}x+y)N_{Q_{t}}(dy),\quad\varphi\in C_b(H),
\;t\geq 0,\;x\in H,
 \end{equation*}
where $N_{Q_t}$ is the Gaussian measure on $H$ of zero mean and covariance operator $Q_t$ (see \cite{DPZ92}).
By Proposition \ref{p.4.2} we know that the OU semigroup $\{R_t\}_{t\geq0}$ is 
a stochastically continuous Markov semigroup in $C_b(H)$.
Moreover, it is well known that for any $t\geq0$, $h\in H$ it holds\footnote{of course, in \eqref{e.35} we consider only the real or the imaginary part}
\begin{equation}\label{e.35}
  R_te^{i\langle \cdot,h\rangle}(x)= e^{i\langle e^{tA}x,h\rangle-\frac12\langle Q_th,h\rangle}, \quad h\in H.
\end{equation}

We denote by $(L,D(L))$ the  infinitesimal generator of $\{R_t\}_{t\geq0}$.
We need the following

\begin{Proposition}  \label{p.1.9}
Let ${\mathcal E}_A(H)$ be the linear span of the real and imaginary part of the functions
$$
   x\mapsto e^{i\langle x,h\rangle},\quad x\in H,\,h \in D(A^*),
$$ 
where $A^*$ is the adjoint of $A$ in $H$. 
For any $\varphi\in C_b(H)$ there exists a three-indexed sequence $(\varphi_{n_1,n_2,n_3})\subset {\mathcal E}_A(H)$ such that
$$
   \lim_{n_1\to\infty}\lim_{n_2\to\infty}\lim_{n_3\to\infty}\varphi_{n_1,n_2,n_3}\stackrel{\pi}{=}\varphi.
$$
Moreover, if $\varphi\in C_b^1(H)$, we have that for any  $h\in H$ it holds
$$
 \lim_{n_1\to\infty}\lim_{n_2\to\infty}\lim_{n_3\to\infty}\langle D\varphi_{n_1,n_2,n_3},h\rangle\stackrel{\pi}{=}\langle D\varphi,h\rangle.
$$
\end{Proposition}
\begin{proof}
Let $\varphi \in C_b(H)$, and let us consider a two-indexed sequence  $(\varphi_{n_1,n_2}) \subset {\cal E}(H)$ as in  Proposition \ref{p.2.1}. 
Let us define the  sequence $(\varphi_{n_1,n_2,n_3})$ by setting
$$
  \varphi_{n_1,n_2,n_3}(x)=\varphi_{n_1,n_2}(n_3R(n_3,A^*)x),\quad x\in H,\,n_3\in\Nset,
$$
where $R(n_3,A^*)$ is the resolvent operator of $A^*$ at $n_3$.
Clearly, $\varphi_{n_1,n_2,n_3}\in {\mathcal E}_A(H)$.
Taking into account that  $nR(n,A^*)x\to x$ as $n\to\infty$ for all $x\in H$, and that for some $c>0$ it holds $|nR(n,A^*)x|\leq c|x|$ for any $x\in H, n\geq 1$, it follows $\varphi_{n_1,n_2,n_3}\stackrel{\pi}{\to}\varphi_{n_1,n_2}$ as $n_3\to \infty$.
If $f\in C^1_b(H)$, we observe that 
$$
\langle D(f(nR(n,A^*)\cdot)(x),h\rangle=\langle Df(nR(n,A^*)x),nR(n,A)h\rangle.
$$
Therefore, be arguing as above, 
we find $\langle D(f(nR(n,A^*)\cdot),h\rangle\stackrel{\pi}{\to} \langle Df(\cdot),h\rangle$ as $n\to\infty$ .
Hence the result follows.
\end{proof}
\begin{Example}
\label{Ex.3.11}
{\em If $A\neq 0$ we have
$
D(L)\cap {\cal E}_A(H)=\{0\}.
$ 
In fact for any $x\in H,h\in D(A^*)$ we have
$$
\lim_{t\to 0^+}\frac{R_t e^{i\langle h,x\rangle}-e^{i\langle h,x\rangle}}{t}=\left[-\frac{1}{2}\;\langle Qh,h\rangle+i
\langle A^*h,x\rangle \right]e^{i\langle h,x\rangle},
$$
which is not bounded when $A\neq 0.$
}
\end{Example}

\begin{Proposition}
\label{p.3.12}
 The set ${\cal I}_A(H)$ 
is $\pi$-dense in $C_b(H)$, it is stable for $R_t$ and ${\cal I}_A(H)\subset D(L)$.
Moreover, it is a $\pi$-core for $(L,D(L))$ and for any $\varphi \in {\cal I}_A(H)$  it holds
\begin{equation} \label{e.36}
L\varphi(x)= \frac12\textrm{Tr}[QD^2\varphi(x)]+ \langle x,A^*D\varphi(x)\rangle,\quad x\in H.
\end{equation}
\end{Proposition}
\begin{proof}
Let  $h\in D(A^*)$ and $a>0$. 
We have  
$$ 
    \lim_{a\to 0^+} \frac1a\int_0^a e^{i\langle
e^{sA}x,h \rangle -\frac12 \langle Q_sh,h\rangle}ds = e^{i\langle
x,h \rangle},\quad x\in H
$$
and
$$
   \sup_{a>0}\bigg|\frac1a\int_0^ae^{i\langle
e^{sA}x,h \rangle -\frac12 \langle Q_sh,h\rangle}ds - e^{i\langle
x,h \rangle}\bigg| \leq 2.
$$
Then ${\cal E}_A(H)\subset \overline{{\cal I}_A(H)}^\pi$.
Consequently, in view of Proposition \ref{p.1.9},  ${\cal I}_A(H)$ is $\pi$-dense in $C_b(H)$.
Now let $t>0$. 
By taking into account \eqref{e.35}, we can apply the Fubini theorem to find
\begin{eqnarray}
 R_t&&\left(\int_0^ae^{i\langle
e^{sA}\cdot,h \rangle -\frac12 \langle Q_sh,h\rangle}ds\right)(x) =\notag\\ 
&&=\int_0^ae^{i\langle
e^{(t+s)A}x,h \rangle -\frac12 \langle Q_te^{sA^*}h,e^{sA^*}h\rangle-\frac12 \langle Q_sh,h\rangle}ds =\notag\\
 &&=\int_0^ae^{i\langle e^{(t+s)A}x,h \rangle -\frac12 \langle Q_{t+s}h,h\rangle}ds=\notag\\
 &&=\int_0^{a+t}e^{i\langle
e^{sA}x,h \rangle -\frac12 \langle Q_sh,h\rangle}ds - \int_0^te^{i\langle
e^{sA}x,h \rangle -\frac12 \langle Q_sh,h\rangle}ds, \label{e.3.13}
\end{eqnarray}
since  $\langle Q_te^{sA^*}h,e^{sA^*}h\rangle = \langle e^{sA}Q_te^{sA^*}h,h\rangle=\langle Q_{t+s}h,h\rangle-\langle Q_sh,h\rangle$. 
Then $R_t ({\cal I}_A(H))$ $\subset {\cal I}_A(H)$.
Now we prove that ${\cal I}_A(H) \subset D(L)$.
Let 
\begin{equation} \label{e.37a}
\DS \varphi(x) =\int_0^ae^{i\langle
e^{sA}x,h \rangle -\frac12 \langle Q_sh,h\rangle}ds.
\end{equation}
By \eqref{e.3.13} we have that
\begin{eqnarray*}
  &&R_t\varphi(x) - \varphi(x) =  \\
  && =\int_a^{a+t}  e^{i\langle
e^{sA}x,h \rangle -\frac12 \langle Q_sh,h\rangle}ds - 
\int_0^te^{i\langle
e^{sA}x,h \rangle -\frac12 \langle Q_sh,h\rangle}ds. 
\end{eqnarray*}
This implies 
\begin{equation} \label{e.39a}
  \lim_{t\to0^+} \frac{R_t\varphi(x) -\varphi(x)}{t} = e^{i\langle
e^{aA}x,h \rangle -\frac12 \langle Q_ah,h\rangle} - e^{i\langle
x,h \rangle} 
\end{equation}
and
$$
   |R_t\varphi(x) -\varphi(x)|\leq 2t.
$$
Then $\varphi \in D(L)$ and  by Proposition \ref{p.4.3} follows that $\mathcal I_A(H)$ is a $\pi$-core for $(L,D(L))$.
In order to prove \eqref{e.36}, it is sufficient take $\varphi$ as in \eqref{e.37a}. 
By a straightforward computation we find that for any $x\in H$ it holds
\begin{eqnarray*}
&&\frac12\textrm{Tr}[QD^2\varphi(x)]+ \langle x,A^*D\varphi(x)\rangle \\
&&\qquad = \int_0^a\left(i\langle A^*e^{sA^*}h,x\rangle-\frac12\langle e^{sA}Qe^{sA^*}h,h \rangle \right)
    e^{i\langle e^{sA}x,h \rangle -\frac12 \langle Q_sh,h\rangle}ds   \\
&&\qquad = \int_0^a \frac{\partial}{\partial s}e^{i\langle e^{sA}x,h \rangle -\frac12 \langle Q_sh,h\rangle}ds\\
&&\qquad = e^{i\langle
e^{aA}x,h \rangle -\frac12 \langle Q_ah,h\rangle} - e^{i\langle
x,h \rangle},
\end{eqnarray*}
cf Example \ref{Ex.3.11}.
By taking into account \eqref{e.39a}, it follows that \eqref{e.36} holds.
\end{proof}

\subsection{Perturbations of the OU operator}
\begin{Proposition} \label{p.5.3}
Under Hypothesis \ref{h.3.1}, let $(L,D(L))$ be the infinitesimal generator of the OU semigroup $\{R_t\}_{t\geq0}$, and let $(K,D(K))$ be the infinitesimal generator of the semigroup $\{P_t\}_{t\geq0}$.
Then  $D(K)\cap C_b^1(H)=D(L)\cap C_b^1(H)$ and for any $\varphi \in D(L)\cap C_b^1(H)$ we have 
  $ K\varphi=L\varphi+\langle D\varphi,F\rangle $.
\end{Proposition}
\begin{proof}
Let $X(t,x)$ be the solution of  equation \eqref{e.4.3} and let us set 
$$
Z_A(t,x)=e^{tA}+\int_0^te^{(t-s)A}Q^{1/2}dW(s).
$$. 
Take $\varphi \in D(L)\cap C_b^1(H)$.
By taking into account that  
$$
   X(t,x)=Z_A(t,x)+ \int_0^te^{(t-s)A}F(X(t,x))ds,
$$
by the Taylor formula we have that $\PP$-a.s. it holds
$$
  \varphi(Z_A(t,x))
   =\varphi(Z_A(t,x))-\varphi(X(t,x))+\varphi(X(t,x))
$$
$$ 
 =\varphi(X(t,x)) -\int_0^1\left\langle D\varphi(\xi Z_A(t,x)+(1-\xi)X(t,x)),\int_0^te^{(t-s)A}F(X(t,x))ds\right\rangle d\xi .
$$
Then we have
$$
   R_t\varphi(x)-\varphi(x) = \EE\big[\varphi(Z_A(t,x))\big]-\varphi(x)= P_t\varphi(x)-\varphi(x)
$$
$$
   -\EE\left[\int_0^1\left\langle D\varphi(\xi Z_A(t,x)+(1-\xi)X(t,x)),\int_0^te^{(t-s)A}F(X(t,x))ds\right\rangle d\xi \right].
$$
Since $\varphi\subset D(L)\cap C_b^1(H)$,
it follows easily that for any $x\in H$
$$
  \lim_{t\to 0^+} \frac{P_t\varphi(x)-\varphi(x)}{t} = L\varphi(x)+\langle D\varphi(x),F(x)\rangle
$$
and 
$$
   \sup_{t\in(0,1]}\bigg\|  \frac{P_t\varphi-\varphi}{t}\bigg\|_0
 \leq \sup_{t\in(0,1]}  \bigg\|  \frac{R_t\varphi-\varphi}{t}\bigg\|_0 +\\
\|D\varphi\|_{C_b(H;\mathcal L(H))}\|F\|_{C_b(H;H)}<\infty,
$$
that implies $\varphi \in D(K)$ and $K\varphi=L\varphi+\langle D\varphi ,F \rangle$.
The opposite inclusion follows by interchanging the role of $R_t$ and $P_t$ in the Taylor formula.
\end{proof}
By the proposition above, we have immediately the following corollary, that proves the first part of Theorem \ref{t.2} 
\begin{Corollary} \label{c.4.6}
Under the hypothesis of Proposition \ref{p.5.3}, we have $\mathcal I_A(H)  \subset D(K)$.
Moreover, the operator $K_0$ is well defined on $\mathcal I_A(H) $ and for any  $\varphi\in \mathcal I_A(H)$ we have $K\varphi=K_0\varphi$.
\end{Corollary}
\begin{proof}
Note that $\mathcal I_A(H) \subset C_b^1(H)$. 
Since by Proposition \ref{p.3.12} we have $\mathcal I_A(H)  \subset D(L)$, by Proposition  \ref{p.5.3} we have $\mathcal I_A(H) \subset D(K)$ and $K\varphi=L\varphi+\langle D\varphi,F\rangle$, for any $\varphi\in \mathcal I_A(H)$.
Finally, by taking into account \eqref{e.36}, it follows that $K\varphi=K_0\varphi$ holds for any $\varphi\in\mathcal I_A(H)$.
\end{proof}
In order to prove that $\mathcal I_A(H)$ is a $\pi$-core for $K$, we need the following approximation result
\begin{Lemma} \label{l.5.4}
Under the hypothesis of Proposition \ref{p.5.3}, let $\varphi\in D(L)\cap C_b^1(H)$.
Then there exists $m\in \Nset$ and an  $m$-indexed sequence $(\varphi_{n_1,\ldots,n_m})\subset \mathcal I_A(H)$ such that
\begin{equation} \label{e.38}
 \lim_{n_1\to\infty} \cdots \lim_{n_m\to\infty}  \varphi_{n_1,\ldots,n_m}\stackrel{\pi}{=}\varphi,
\end{equation}
\begin{equation}\label{e.39}
 \lim_{n_1\to\infty} \cdots \lim_{n_m\to\infty} \frac12\textrm{Tr}\big[QD^2\varphi_{n_1,\ldots,n_m}\big]+\langle \cdot,A^*D\varphi_{n_1,\ldots,n_m}\rangle  \stackrel{\pi}{=} L\varphi,
\end{equation}
and for any $h\in H$
\begin{equation} \label{e.40}
  \lim_{n_1\to\infty} \cdots \lim_{n_m\to\infty}  \langle D \varphi_{n_1,\ldots,n_m},h\rangle \stackrel{\pi}{=}
    \langle D  \varphi,h\rangle.
\end{equation}
\end{Lemma}
\begin{proof}
We observe that the results of Proposition \ref{p.1.9} holds also by approximations with functions in $\mathcal I_A(H)$.
Indeed, let $(\varphi_{n_1,n_2,n_3 })\subset \mathcal E_A(H)$ as in Proposition \ref{p.1.9}.
By setting, for any $n_1, n_2,n_3,n_4\in\Nset$
$$
  \varphi_{n_1,n_2,n_3,n_4 }(x)=n_4\int_0^\frac{1}{n_4}R_t\varphi_{n_1,n_2,n_3 }(x)dt
$$
we have, thanks to \eqref{e.35}, that $\varphi_{n_1,n_2,n_3,n_4 }\in \mathcal I_A(H)$.
Clearly, 
$$
\lim_{n_1\to\infty}\cdots \lim_{n_4\to\infty}\varphi_{n_1,n_2,n_3,n_4 }\stackrel{\pi}{=}\varphi.
$$
Moreover, since $D (R_t f)=e^{tA^*}R_t(D\varphi)$ (cf, e.g., \cite[Proposition 6.2.9]{DPZ02}), we find that for any $h\in H$ it holds
$$
   \langle D \varphi_{n_1,n_2,n_3,n_4 }(x),h\rangle=
   n_4\int_0^\frac{1}{n_4}R_t\big(\langle D\varphi_{n_1,n_2,n_3 }(\cdot),e^{tA}h\rangle\big)(x) dt.
$$
Hence, 
$$
\lim_{n_1\to\infty}\cdots \lim_{n_4\to\infty} \langle D\varphi_{n_1,n_2,n_3,n_4 },h\rangle\stackrel{\pi}{=}\langle D\varphi,h\rangle.
$$

Now we construct the desired approximation for $\varphi\in D(L)\cap C_b^1(H)$.
Let $\varphi\in D(L)\cap C_b^1(H)$ and $(\varphi_{n_2})\subset \mathcal I_A(H)$ as above 
(of course, for simplicity we assume that this approximation has only one index; this does not reduce the generality of the proof).
By setting $(\varphi_{n_1,n_2,n_3})$ as in \eqref{e.32a} with $R_t$ instead of $P_t$, we have that \eqref{e.38}, \eqref{e.39} hold, by the same argument of the proof of Proposition \ref{p.4.3}.

We now observe that for any $n_1,n_2,n_3\in \Nset$, the function $\varphi_{n_1,n_2,n_3}$ is differentiable in every $x\in H$ along any direction $h\in H$, with differential
$$
   \langle D \varphi_{n_1,n_2,n_3}(x),h\rangle=  
   \varphi_{n_1,n_2,n_3}(x)=\frac{1}{n_3}\sum_{i=1}^{n_3} R_{\frac{i}{n_1n_3}}\big(\langle D\varphi_{n_2}(\cdot),e^{\frac{i}{n_1n_3}A }h\rangle\big)(x)
$$
Moreover,
$$
 \sup_{n_1,n_2,n_3\in \Nset }\| \langle D \varphi_{n_1,n_2,n_3},h\rangle\|_0\leq 
  \sup_{n_2} \| D\varphi_{n_2}\|_{C_b(H;H)}\sup_{0\leq t\leq 1}\|e^{tA}\|_{\mathcal L(H)}|h| <\infty.
$$
Now by arguing as for Proposition  \ref{p.4.3}, it yields \eqref{e.40}.
\end{proof}

\subsection{The case $F\in C_{\lowercase{b}}^2(H;H)$}
The following proposition is proved in \cite[section 3.3]{DP04}.
\begin{Proposition} \label{p.4.8}
Let us assume Hypothesis \ref{h.3.1} and that $F\in C_{\lowercase{b}}^2(H;H)$, that is $F:H\to H$ is two time differentiable with bounded differentials.
Then the semigroup $\{P_t\}_{t\geq0}$ defined in \eqref{e.38a} maps $C_b^1(H)$ into $C_b^1(H)$, and for any $f\in C_b^1(H)$, $h\in H$ we have
$$
    \langle DP_tf(x),h\rangle=\EE\big[ \langle Df(X(t,x)),\eta^h(t,x)\rangle  \big],
$$
where $\eta^h(t,x)$ is the mild solution of the differential equation in $H$
$$
  \begin{cases} 
\eta^h(t,x)= A\eta^h(t,x)+\langle DF(X(t,x)),\eta^h(t,x)\rangle,& t>0, \\
   \eta^h(0,x)=h. &{}
  \end{cases}
$$
\end{Proposition}
\begin{Corollary} \label{c.4.12}
Under the hypothesis of Proposition \ref{p.4.8}, let $(K,D(K))$ be the infinitesimal generator of  $\{P_t\}_{t\geq0}$.
Then, for any $\lambda>0,\omega+M\|DF\|_0$, the resolvent $R(\lambda,K)$ of $K$ at $\lambda$ maps $C_b^1$ into $C_b^1(H)$ and it holds
\begin{equation} \label{e.43a}
  \|DR(\lambda,K)f\|_{C_b(H;H)}\leq \frac{M\|Df\|_{C_b(H;H)}}{\lambda-(\omega+M\|DF\|_{C_b(H;\mathcal L(H))})},\quad f\in C_b^1(H).
\end{equation}
\end{Corollary}
\begin{proof}
Let  $f\in C_b^1(H)$.
For any $t\geq0$, $P_tf\in C_b^1(H)$ and for any $x,h\in H$ it holds 
$$
   \langle DP_tf(x),h \rangle =\EE\big[ \langle Df(X(t,x)),\eta^h(t,x)\rangle\big],
$$
where  $\eta^h(t,x)$ is as in Proposition \ref{p.4.8}.
It is also easy to see that\footnote{in order to avoid heavy notations we set $\|DF\|=\|DF\|_{C_b(H;\mathcal L(H))}$} 
$$
    |\eta^h(t,x)|\leq Me^{(\omega+M\|DF\|)t}|h|,
$$
see, e.g., \cite[Theorem 3.6]{DP04}.
Hence, by (vi) of Theorem \ref{t.2.9}, we have  
\begin{eqnarray*}
   |\langle DR(\lambda,K)f(x),h\rangle |&=&\left|\int_0^\infty e^{-\lambda t} \EE\big[ \langle Df(X(t,x)),\eta^h(t,x)\rangle\big]dt \right| \\
   &\leq& M  \|Df\|_{C_b(H;H)} \int_0^\infty e^{-\lambda t} e^{(\omega+M\|DF\|)t}|h|dt \\
&=&\frac{M\|Df\|_{C_b(H;H)}}{\lambda-(\omega+M\|DF\|_{C_b(H;\mathcal L(H))})}|h|,
\end{eqnarray*}

for any $h \in H$. 
Therefore,  \eqref{e.43a}  follows.
\end{proof}

\begin{Proposition} \label{p.6.2}
 Let us assume that that Hypothesis \ref{h.3.1} hold and let $F\in C_b^2(H;H)$.
Denoted by $\{P_t\}_{t\geq0}$ the transition semigroup defined in \eqref{e.38a}, let $(K,D(K))$ be its infinitesimal generator.
Then, the set $\mathcal I_A(H)$ introduced in Theorem \ref{t.1.4} is a $\pi$-core for $(K,D(K))$, and for any $\varphi\in D(K)$ there exists $m\in \Nset$ and an $m$-indexed sequence $( \varphi_{n_1,\ldots,n_m})\subset \mathcal I_A(H)$ such that
\begin{equation}\label{e.43}
 \lim_{n_1\to\infty} \cdots \lim_{n_m\to\infty} K_0\varphi_{n_1,\ldots,n_m}\stackrel{\pi}{=} K\varphi.
\end{equation}
\end{Proposition}
\begin{proof} 
Let $\varphi\in D(L)\cap C_b^1(H)$.
By Proposition \ref{p.5.3} we have that $\varphi\in D(K)\cap  C_b^1(H)$.
Hence, by (i) of Theorem \ref{t.2.9} we have $P_t\varphi\in D(K)$ and by Proposition \ref{p.4.8} we have $P_t\varphi\in  C_b^1(H)$, for any $t\geq0$. 
So  $P_t\colon  D(L)\cap C_b^1(H)$ $\to D(L)\cap C_b^1(H)$, for any $t\geq0$.
Moreover, $\mathcal I_A(H)\subset D(L)\cap C_b^1(H)$ and so $D(L)\cap C_b^1(H)$ is $\pi$-dense in $C_b(H)$, in view of the fact that $\mathcal I_A(H)$  is  $\pi$-dense in $C_b(H)$ (cf Prop. \ref{p.3.12}).
Therefore, by Proposition \ref{p.4.3}, $D(L)\cap C_b^1(H)$ is a $\pi$-core for $(K,D(K))$.
So there exists a sequence $(\varphi_m)\subset \mathcal I_A(H)$ (we assume that the sequence has one index) such that $L\varphi_m+\langle D\varphi_m,F\rangle\stackrel{\pi}{\to}K\varphi$, as $m\to \infty$.
Now, thanks to Lemma \ref{l.5.4}, we can approximate any $\varphi_m$ by a sequence $(\varphi_{m,n}) \subset \mathcal I_A(H)$
in such a way that $\varphi_{m,n}\stackrel{\pi}{\to}\varphi_m$, $L\varphi_{m,n}\stackrel{\pi}{\to}L\varphi_m$ as $n\to\infty$ and
$ \langle D \varphi_{m,n},h\rangle \stackrel{\pi}{\to} \langle D \varphi_m,h\rangle$ as $n\to\infty$, for any $h\in H$.
Since $F:H\to H$ is bounded, we have $ \langle D \varphi_{m,n},F\rangle \stackrel{\pi}{\to} \langle D \varphi_m,F\rangle$ as $n\to\infty$.
Finally, since $\varphi_{m,n}\in \mathcal I_A(H) $ by Corollary \ref{c.4.6} it follows  \eqref{e.43}.
\end{proof}
\subsection{The Lipschitz case and conclusion of the proof}
%
%
Corollary \ref{c.4.6} proves that $K$ is an extension of $K_0$, and that $K\varphi=K_0\varphi$, $\forall \varphi \in \mathcal I_A(H)$.
It remains to prove that $\mathcal I_A(H)$ is a $\pi$-core for $K$.

We  denote by $L_F$ the Lipschitz constant of $F$.
Let $\varphi \in D(K)$, $\lambda > \max\{0,\omega+L_F\}$ and set $f=\lambda\varphi-K\varphi$.
Since $C_b^1(H)$ is dense in $C_b(H)$ with respect to the supremum norm (see \cite{LL86}), there exists a sequence $(f_{n_1})\subset C_b^1(H)$ such that $\|f_{n_1}-f\|_0\to0$ as $n_1\to \infty$.
Clearly, if $\varphi_{n_1}=R(\lambda,K)f_{n_1}$ we have
\begin{equation} \label{e.44b}
    \lim_{n_1\to\infty}K\varphi_{n_1}\stackrel{\pi}{=}K\varphi.
\end{equation} 
Now we consider a sequence of functions $(F_{n_2})_{n_2\in\Nset}\subset C_b^2(H;H)$ such that 
\begin{equation} \label{e.46a}
    \lim_{n_2\to\infty}F_{n_2}(x)=F(x),\quad \forall x\in H
\end{equation}
and
\begin{equation}\label{e.47a}
 \sup_{{n_2}\in\Nset}\|F_{n_2}\|_{C_b(H;H)}\leq \|F\|_{C_b(H;H)},\quad \sup_{n_2\in\Nset}  \|DF_{n_2}\|_{C_b(H;\mathcal L(H))}\leq L_F.
\end{equation}
This construction is not too difficult but technical and an example can be found in \cite[section 3.3.1]{DP04}.
Let $X^{n_2}(t,x)$ be the solution  of \eqref{e.4.3} with $F_{n_2}$ instead of $F$. 
It is straightforward to see that for any $T>0$, $x\in H$
$$
   \lim_{{n_2}\to\infty} \sup_{t\in [0,T]}\EE\big[|X^{n_2}(t,x)-X(t,x)|^2  \big] = 0.
$$ 
Hence, if $P_t^{n_2}$ is the transition semigroup associated to $X^{n_2}(t,x)$, we have that for any $\varphi\in C_b(H)$
$$
   \lim_{{n_2}\to\infty} P_t^{n_2}\varphi\stackrel{\pi}{=}P_t\varphi.
$$
We denote by $(K_{n_2},D(K_{n_2}))$ the infinitesimal generator of the transition semigroup $\{P_t^{n_2}\}_{t\geq0}$, as in \eqref{e.3}.
We also set
$$
   K_{0, n_2}\varphi(x)=K_0\varphi(x)+\langle D\varphi(x),F_{n_2}-F(x)\rangle,\,\varphi\in \mathcal I_A(H),\,x\in H.
$$
If $R(\lambda,K_{n_2})$ is the resolvent of $K_{n_2}$ at $\lambda$ (cf (vi) of Theorem \ref{t.2.9}), we have 
\begin{equation*} 
     \lim_{n_2\to\infty}R(\lambda,K_{n_2})f\stackrel{\pi}{=}R(\lambda,K)f,
\end{equation*}
for any $f\in C_b(H)$.
Setting $\varphi_{n_1,n_2}=R(\lambda,K_{n_2})f_{n_1}$, for any $n_1\in\Nset$ we have
\begin{equation} \label{e.44a}
\lim_{n_2\to\infty}\varphi_{n_1,n_2}\stackrel{\pi}{=} \varphi_{n_1},\quad \lim_{n_2\to\infty}K_{n_2}\varphi_{n_1,n_2}\stackrel{\pi}{=} K\varphi_{n_1}.
\end{equation}
Moreover, since $F_{n_2}\in C_b^2(H;H)$, by Corollary \ref{c.4.12} we have that $R(\lambda,K_{n_2}):C_b^1(H)\to C_b^1(H)$ and
$$
  \|D\varphi_{n_1,n_2} \|_{C_b(H; H)} \leq   \frac{M\| D\varphi_{n_1} \|_{C_b(H; H )}}{\lambda-(\omega+\|DF_{n_2}\|_{C_b(H;\mathcal L(H))} )}\leq   \frac{M\| D\varphi_{n_1} \|_{C_b(H;H)}}{\lambda-(\omega+L_F )},
$$
for any $n_1, n_2\in \Nset$. 
Consequently, by \eqref{e.46a}, \eqref{e.47a} it follows
\begin{equation} \label{e.45}
   \lim_{n_2\to\infty}\langle D\varphi_{n_1,n_2},F-F_{n_2}\rangle\stackrel{\pi}{=}0.
\end{equation}
Since $f_{n_1}\in C_b^1(H)$, by Corollary \ref{c.4.12} we have $\varphi_{n_1,n_2}\in D(K_{n_2})\cap C_b^1(H)$.
By  Proposition \ref{p.6.2}, for any $n_1$, $n_2\in \Nset$ we can find a sequence $(\varphi_{n_1,n_2,n_3})\subset \mathcal I_A(H)$ such that
\begin{equation} \label{e.47}
   \lim_{n_3\to \infty} K_{0,n_2}\varphi_{n_1,n_2,n_3}\stackrel{\pi}{=}L\varphi_{n_1,n_2}+\langle D\varphi_{n_1,n_2},F_{n_2}\rangle=K_{n_2}\varphi_{n_1,n_2}.
\end{equation}
Hence we have
$$
  K_0\varphi_{n_1,n_2,n_3}=K_{0,n_2}\varphi_{n_1,n_2,n_3}+\langle D\varphi_{n_1,n_2,n_3},F-F_{n_2}\rangle
$$
and by \eqref{e.44a}, \eqref{e.45}, \eqref{e.47} it follows
$$
\lim_{n_2\to \infty}\lim_{n_3\to \infty} K_0\varphi_{n_1,n_2,n_3}
$$
$$
  \stackrel{\pi}{=}
  \lim_{n_2\to \infty}K_{n_2}\varphi_{n_1,n_2}+\langle D\varphi_{n_1,n_2},F-F_{n_2}\rangle  \stackrel{\pi}{=}K \varphi_{n_1}.
$$
Now the result follows by \eqref{e.44b}. \qed
%
%
\section{Proof of Theorem \ref{t.1.4} }
Let $\mu\in \mathcal M(H)$ and assume that $\{\mu_t\}_{t\geq0}$ is a solution of the measure equation \eqref{e.2.5a}.
Denoting by $(K,D(K))$  the infinitesimal generator of the semigroup \eqref{e.38a},
by Theorem \ref{t.2} we have that $\mathcal I_A(H)$ is a $\pi$-core for $(K,D(K))$, and that $K\varphi=K_0\varphi$, for any $\varphi \in \mathcal I_A(H)$.
This implies that $(K,D(K))$ is an extension of $K_0$.
So it is easy to see that $\{P_t^*\mu\}_{t\geq0}$ is a solution of the measure equation \eqref{e.2.5a}.
Hence, if $\varphi\in D(K)$ there exists a sequence\footnote{For simplicity we assume that this sequence has only one index} $(\varphi_n)\subset \mathcal I_A(H)$ such that
$$
  \lim_{n\to\infty}\varphi_n\stackrel{\pi}{=}\varphi,\quad \lim_{n\to\infty}K_0\varphi_n\stackrel{\pi}{=}K\varphi.
$$ 
For any $t \geq0$ we find
\begin{eqnarray*}
   \int_H\varphi(x)\mu_t(dx)-\int_H\varphi(x)\mu (dx)&=&\lim_{n\to\infty}
   \left( \int_H\varphi_n(x)\mu_t(dx)-\int_H\varphi_n(x)\mu (dx)\right)\\
  &=&  \lim_{n\to\infty}\int_0^t\left( \int_HK_0\varphi_n(x)\mu_s(dx)\right)ds.
\end{eqnarray*}
Now observe that for any $s\geq0$ it holds
$$
\lim_{n\to\infty}\int_HK_0\varphi_n(x)\mu_s(dx)=\int_HK\varphi(x)\mu_s(dx)
$$
and 
$$
\left|\int_HK_0\varphi_n(x)\mu_s(dx)\right|\leq \sup_{n\in\Nset}\|K_0\varphi_n\|_0\|\mu_s\|_{TV}.
$$
Hence, by taking into account \eqref{e.5b} and that $\sup_{n\in\Nset}\|K_0\varphi_n\|_0<\infty$, we can apply the dominated convergence theorem to obtain 
$$
    \lim_{n\to\infty}\int_0^t\left( \int_HK_0\varphi_n(x)\mu_s(dx)\right)ds=
   \int_0^t\left( \int_HK\varphi(x)\mu_s(dx)\right)ds
$$
So, $\{\mu_t\}_{t\geq0}$ is also a solution of the measure equation for $(K,D(K))$.
Since by Theorem \ref{t.1} such a solution is unique, if follows that  the measure equation \eqref{e.2.5a} has a unique solution, done by $\{P_t^*\mu\}_{t\geq0}$. \qed
\begin{Remark}{\em
If  $(K,D(K))$ is the infinitesimal generators of a stochastically continuous Markov semigroup and  
 $D$ is a $\pi$-core for $(K,D(K))$, we can extend the theorem above to the operator   $K_0:=K|_D$.
Indeed, all the computations are similar.
}
\end{Remark}

\bibliographystyle{amsplain}
\bibliography{biblioluigi}
\end{document}